\documentclass[a4paper, 12pt]{elsarticle}

\usepackage{ucs}
\usepackage[utf8x]{inputenc}
\usepackage{amsmath,amssymb,amsfonts,mathtools}
\usepackage[english]{babel}
\usepackage{fontenc}
\usepackage{graphicx}
\usepackage[caption=false]{subfig}
\usepackage{ulem} 
\usepackage[margin=0pt,font=footnotesize, labelsep=colon]{caption}
\usepackage{tabularx}
\usepackage{multirow}
\usepackage{siunitx}
\usepackage{booktabs}
\usepackage{color}
\usepackage[table]{xcolor}
\usepackage{colortbl}
\usepackage{supertabular}
\usepackage{array}

\definecolor{rgrey}{gray}{0.75}
\newcommand{\mean}[1]{\langle #1 \rangle}

\begin{document}

\begin{frontmatter}

\title{Transmission needs across a fully renewable European power system}
  
\author[label1]{Rolando A. Rodriguez}
\ead{rar@imf.au.dk}
\author[label2]{Sarah Becker}
\author[label1]{Gorm B. Andresen}
\author[label3]{Dominik Heide}
\author[label1]{Martin Greiner}
\address[label1]{Department of Engineering and Department of Mathematics, Aarhus University, Ny Munkegade 118, 8000 Aarhus C,  Denmark}
\address[label2]{Frankfurt Institute for Advanced Studies (FIAS), Johann Wolfgang Goethe Universit{\"a}t, Ruth-Moufang-Stra{\ss}e 1, 60438 Frankfurt am Main, Germany}
\address[label3]{Deutsches Zentrum f\"ur Luft- und Raumfahrt (DLR), Pfaffenwaldring 38-40, 70569 Stuttgart, Germany}
  
\date{\today}

\begin{abstract}
The residual load and excess power generation of 27 European countries with a 100\% penetration of variable renewable energy sources are explored in order to quantify the benefit of power transmission between countries. Estimates are based on extensive weather data, which allows for modelling of hourly mismatches between the demand and renewable generation from wind and solar photovoltaics. For separated countries, balancing is required to cover around 24\% of the total annual energy consumption. This number can be reduced down to 15\% once all countries are networked together with unconstrained interconnectors. The reduction represents the maximum possible benefit of transmission for the countries. The total Net Transfer Capacity of the unconstrained interconnectors is roughly twelve times larger than current values. However, constrained interconnector capacities six times larger than the current values are found to provide 97\% of the maximum possible benefit of cooperation. This motivates a detailed investigation of several constrained transmission capacity layouts to determine the export and import capabilities of countries participating in a fully renewable European electricity system.
\end{abstract}


\end{frontmatter}
\section{Introduction}
\label{sec:one}

The sustainability of the world's energy supply is strongly dependent on the successful integration of renewable sources. Variable Renewable Energy Sources (VRES), such as wind and solar energy, promise to be key elements in future energy systems \cite{Jacobson:2011a, Jacobson:2011b, epia, Greenpeace, Cho2010}. The nature of VRES makes them hard to integrate into an electrical system that was built on more or less predictable loads with dispatchable generation. In small penetrations, the variations can be absorbed without much consequence, but will be harder to ignore in a future, highly renewable, macro energy system. The spatio-temporal dispersion of the weather patterns that define the output of wind and solar energy will lead to fluctuating mismatches between regional demand for and generation of electricity. This will give rise to new challenges for countries with a high penetration of VRES, such as the need for back-up conventional balancing, flexible demand, dispatchable renewable sources, such as hydroelectric reservoirs or biomass, increased transmission capacities to neighboring regions and energy storage \cite{Jacobson:2011b, Greenpeace}. In order to understand and to design the future energy systems with dominant shares of VRES, we need to let the weather decide.

For the optimal integration of VRES in future 100\% renewable electricity systems, one wishes to make as much use as possible of renewables while minimizing the need for conventional balancing, both in the installed power capacity required and the energy expended \cite{gormconf}. Additionally, we wish to minimize the need for storage \cite{Heide:2010fk, Heide:2011kx} and transmission capacities \cite{Schaber, Schaber2}. In determining lower bounds on the need for storage and transmission, the synergies between these factors and the need for balancing must be well understood \cite{Rasmussen2012}. In this article, we focus on determining the synergy between transmission and balancing.

Several studies have assessed the need for a larger transmission network \cite{ Greenpeace, Tradewind, Buijs, energynautics}. Despite the planned investments in grid strength, the European Network of Transmission System Operators for Electricity (ENTSO-E) has identified 100 bottlenecks in their network development plan \cite{ENTSOEtyndp}, with 80\% of them due to integration of renewables. By looking at characteristic weather patterns and possible wind and solar power generation across Europe, potential transmission between regions have been estimated. This has been done for Germany \cite{Zugno}, and with an economic approach for Europe \cite{Schaber, Steinke}. A similar study has looked at regional aggregation and transmission in the United States \cite{Corcoran}. Estimates on the size of an ideal transmission grid are large, such as 20 GW for the link between France and Spain \cite{energynautics}, which is over 15 times larger than the current interconnector capacity. There is a conflict between the need for maximizing the integration of fluctuating VRES and minimizing the expansion of the transmission system. 

Starting from the same large weather database as presented in \cite{Heide:2010fk, Heide:2011kx, Rasmussen2012}, we estimate the potential output of wind and solar photovoltaic energy for any given country in a 27-node representation of Europe. In Section 2, we introduce a model which calculates the local mismatches between VRES generation and load in this set of interconnected countries, and which distributes the excess generation in a way that maximizes the use of renewables. An efficient usage of renewables also minimizes the need for balancing energy $E_B$ coming from conventional dispatchable resources. Section 2 also explains how this optimal distribution of VRES excess generation can be found by performing a novel DC power flow calculation with constrained interconnectors. It also determines the interplay between installed transmission capacity and the benefit coming from transmission. In Section 3, the DC power flow model is applied to the case of a future, 100\% renewable Europe. A minimum $E_B$ that each country can attain through an optimal mix of wind and solar is found, and then compared to that of a fully connected, unconstrained Europe. The total $E_B$ resulting from this unconstrained flow leads to the maximum benefit of transmission, when countries can make the most use of the renewable excess generation of their neighbors. By applying the constrained DC power flow calculation we find a precise relation between the installed transmission capacity and the required total balancing energy $E_B$. Section 4 discusses the limits to import and export capabilities, and the reduction of conventional power capacities. The conclusion is presented in Section 5.

\section{Methodology}
\label{sec:two}

The following is a method to determine power flows in a power system with a large amount of VRES generation, and the benefit they bring by reducing the need for balancing. Power flow calculations are detailed for unconstrained and constrained cases.

\subsection{Definitions}

For a node $n$, representing a country, the hourly VRES generation and the electrical load will generally not be equal. The hourly mismatch between the load $L_n$ and the combined output of wind $G^W_n$ and solar $G^S_n$ generation is defined as \cite{Heide:2010fk}
\begin{equation}
\label{eq:two1}
  \Delta_n(t)
    =  \gamma_n \left ( 
           \alpha^W_n \frac{G^W_n(t)}{\mean{G^W_n(t)}} + (1-\alpha^W_n) \frac{G^S_n(t)}{\mean{G^S_n (t)}}
        \right ) \cdot \mean{L_n(t)} 
        - L_n(t)        \; .
\end{equation}
Here, $t$ represents the hourly timestep, $\gamma_n$ the total VRES penetration and $\alpha^W_n$ the relative mix between wind and solar generation at node $n$. The VRES generation is normalized to the mean value of the load. For the case when $\gamma_n=1$, VRES generate as much energy, on average, as is consumed by the load. The mean mismatch is zero, but, due to the fluctuations of the generation and the load, the mismatch will almost always be either positive in case of excess generation or negative in case of deficit generation. 

The negative part of the mismatch defines the residual load of a country,
\begin{equation}
\label{eq:two2}
  \Delta^-_n(t)
    =  \max\left\{ -\Delta_n(t) , 0 \right\}
        \; , 
\end{equation}
which needs to be balanced by additional power generation sources other than wind and solar. The positive part of the mismatch is excess power  
\begin{equation}
\label{eq:two3}
  \Delta^+_n(t)
    =  \max\left\{ \Delta_n(t) , 0 \right\}
        \; ,
\end{equation}
which can be exported or needs to be curtailed. The time averages of (2) and (3) are denoted $\langle \Delta_n^- \rangle$ and $\langle \Delta_n^+ \rangle$. For $\gamma_n=1$, these two averages are identical. 

\subsection{Unconstrained DC power flow}

We now assume that the nodes are joined by links that allow the transfer of energy from one to another. Given a directed graph consisting of $N$ nodes and $L$ links, the topology can be described by the $N \times L$ incidence matrix
\begin{equation}
\label{eq:incidence}
  K_{n,l}
    =  \left \{ \begin{array}{rl}
        1  & \text{if link} \ l \ \text{starts at node} \ n, \\
        -1 & \text{if link} \ l \ \text{ends at node} \ n, \\
        0  & \text{otherwise} \: .
        \end{array} \right. 
\end{equation}
The choice of link direction is arbitrary, but once taken it has to be kept fixed throughout the calculations.
To calculate the occurring flows and the balancing needs, we employ the DC approximation to the full AC power flow \cite{WoodWollenberg1996}. It is valid as long as the network is in steady state, the resistances of the links can be neglected and no significant voltage phase shifts occur between the nodes \cite{Hertem}. For convenience, all susceptances are assumed to be uniform and equal to one.

When the combined mismatch of all nodes is exactly zero, the total excess power matches the total residual load in the network. In this case, the unconstrained DC power flows $F_l$ follow from the equation
\begin{equation}
 \Delta_n(t)-\sum_{l=1}^L K_{n,l} F_l(t) = 0 
\end{equation}
and the minimization of the square sum $\sum_{l=1}^L F_l^2$ \cite{WoodWollenberg1996, Hertem}.

In general, the combined mismatch of all nodes is not zero. In this case, the balancing at each node can be expressed as the negative part of the local mismatch $\Delta_n$ minus the net exports $\sum_{l=1}^L K_{n,l}F_l$:
\begin{equation}
\label{eq:Bn}
  B_n(t)  =  - \min\left\{
                  \left[ \Delta_n(t) - \sum_{l=1}^L K_{n,l} F_l(t) \right]
                  \; , \; 0 \;
                  \right\}
                  \; .
\end{equation}
Likewise, the excess power, which has to be curtailed or used otherwise, can be expressed as the positive part of the difference between the mismatch and the net exports:
\begin{equation}
  C_n(t)  =  \max\left\{
                  \left[ \Delta_n(t) - \sum_{l=1}^L K_{n,l} F_l(t) \right]
                  \; , \; 0 \;
                  \right \}
                  \; .
\end{equation}

The balancing $B_n$, the excess power $C_n$ and the flows $F_l$ at all nodes and links are determined by a two-step optimization procedure. The first priority is the minimization of the overall balancing for each hour:
\begin{equation}
\label{eq:three2}
  B_{\rm min}(t)
     =   \min_{F_l} \sum_{n=1}^N B_n(t)
          \; ,
\end{equation}
which guarantees an optimal usage of VRES across all nodes. Since this does not yet determine the flows in a unique manner, they are fixed in a second step,
\begin{equation}
\label{eq:uncprob}
\begin{aligned}
&\min_{F_l} &&  \sum_{l=1}^L F_l^2 \\
&\text{s.t.}   &&  \sum_{n=1}^N B_n = B_{\rm min} \; ,
\end{aligned}
\end{equation}
which minimizes the quadratic flows with the constraint of keeping the total balancing at its minimal value found in the first step. The two steps ensure that we arrive at the most localized minimal flows which allow an optimal sharing of renewables.

\subsection{Constrained DC power flow}

Today's power grids are constrained by the Net Transfer Capacities (NTC). These constraints are based not only on the physical properties of the interconnectors, the Total Transfer Capacities (TTC), but also on the strength of the grid on either sides of the link and on the security policies of the participating countries \cite{SCOPF}. As a result, transmission over links is constrained with different values in each direction.

The DC power flow as stated in equation (\ref{eq:uncprob}) is easily constrained by adding limits $f_l^- \leq F_l \leq f_l^+$:
\begin{equation}
\label{eq:conprob}
\begin{aligned}
  \text{Step 1: }&\qquad&\min_{F_l} \quad & \sum_{n=1}^N B_n  \equiv B_{\rm min} &\\
   &\qquad& \text{s.t.} \quad & f_l^{-}\leq F_l \leq f_l^+ & \\
  \text{Step 2: }&\qquad&\min_{F_l} \quad& \sum_{l=1}^L F_l^2 & \\
   &\qquad& \text{s.t.} \quad& f_l^{-}\leq F_l \leq f_l^+ & \\
   & \qquad&               &  \sum_{n=1}^N  B_n = B_{\rm min}  \; .  & \\
\end{aligned}
\end{equation}
The first step can be solved linearly with the help of a slack variable, whereas the second step is a quadratic programming optimisation problem. Solutions can be obtained with several computational solving tools. The high speed of the solvers generated by CVXGEN \cite{CVXGEN} made this tool well suited for our purposes.

\subsection{Benefit of transmission}

The constrained optimisation (\ref{eq:conprob}) determines the power flows $F_l$ on all links and, via (\ref{eq:Bn}), also the residual loads $B_n$ on all nodes. It is important to note that the resulting power flows and residual loads depend on the transmission capacity layout $\{ f_l^\pm \}$, which constrains the power flows. A measure on how much a transmission capacity layout is able to reduce the overall residual load is defined by the benefit of transmission
\begin{equation}
\label{eq:benefit}
  \beta  =  \frac{ E_B(0) - E_B(\{f_l^\pm\}) }{ E_B(0) - E_B(\infty) }
                \; .
\end{equation} 
The average annual balancing energy 
\begin{equation}
  E_B  =  N_h \sum_{n=1}^N \langle B_n \rangle
\end{equation}
is defined as the sum over the country-specific average residual loads, multiplied with the total number $N_h=8760$ of hours per year. The case of zero transmission capacity layout $f^\pm=0$ results in zero power flows, so that the average overall annual balancing energy 
\begin{equation}
  E_B(f_l^\pm\!=\!0)  =  N_h \sum_{n=1}^N \langle \Delta_n^- \rangle
  \end{equation} 
can be expressed as the country sum over the individual negative mismatches (\ref{eq:two2}). Its benefit of transmission is $\beta=0$. Unconstrained power flows represent the other extreme and are the result of an infinitely strong layout. The resulting annual balancing energy 
\begin{equation}
  E_B(\infty)  =  N_h \left\langle \max\left( -\sum_{n=1}^N \Delta_n \; , \; 0 \; \right) \right\rangle
\end{equation}
can be obtained without a flow calculation, and leads to the maximum benefit of transmission $\beta=1$. Each transmission capacity layout $\{f_l^\pm\}$ between the two extremes will result in an annual balancing energy $E_B(\infty) \leq E_B(f_l^\pm) \leq E_B(0)$ and a benefit of transmission $0\leq\beta\leq1$.

\section{Case study: Power transmission in a fully renewable Europe}
\label{sec:three}

The methodology developed in the previous section is applied to the case of a 27-node fully renewable European network. Each country has a penetration $\gamma_n = 1$ of combined wind and solar power generation, which on average is equal to its load. Country-specific optimal mixes are first determined. Then, the unconstrained formulation (\ref{eq:uncprob}) for the power flow is used to estimate maximum transmission capacities on the interconnectors. Finally, the constrained power flow problem (\ref{eq:conprob}) is solved with different capacity layouts and the respective benefits of transmission $\beta$ are determined.

\subsection{Country specific optimal mixes}

Focusing on wind and solar energy, VRES generation potentials were determined for 27 European countries from a large weather database, spanning eight years from 2000 to 2007 \cite{Heide:2010fk}. This includes both on- and off-shore regions with a spatial resolution of $47$ km $\times$ $47$ km. Weather measurements were used to determine wind and solar energy generation time series with an hourly resolution. To obtain the absolute power output for a country, capacity scaling factors were applied to each grid cell belonging to the country to be aggregated. The scaling factors reflect the assumed installed wind and solar capacity at that location. 
Historical data for the electricity demand was used to generate hourly load time series for all 27 countries covering the same 8 years from 2000 to 2007. The time series for each country were then de-trended to correct for the approximately 2\% annual increase in  electricity demand. Average loads for the countries can be seen in Table 1. 

\begin{table}
	\label{tab:1}
	\centering
    	\scriptsize
    \begin{tabularx}{135mm}{  llS  llS  llS }
    \toprule
    ISO&Country & {$\langle L \rangle$} & ISO&Country & {$\langle L \rangle$}  &ISO&Country & {$\langle L \rangle$} \\ 
     && {(GW)} & && {(GW)} &&& {(GW)}\\ [0.50 ex] 
    \midrule
        	DE &	 Germany & 54.2	& FI & Finland & 9.0 & RS 	& Serbia & 3.9 \\	
	FR &	 France 	& 51.1	& CZ & Czech Republic & 6.7	& IE & Ireland & 3.2 \\	
	GB &	 Great Britain 	& 38.5 & AT 	& Austria & 5.8 & SK & Slovakia & 3.1 \\	
	IT &	 Italy 	& 34.5	& GR & Greece & 5.8 & BA & Bosnia \&	 Herz.  &  3.1 \\
	ES &	 Spain 	& 24.3	& RO & Romania & 5.4	 & HR & Croatia & 1.6 \\	
	SE &	 Sweden 	& 16.6	& BG & Bulgaria & 5.1 & SI & Slovenia & 1.4 \\	
	PL &	 Poland 	& 15.2	& PO & Portugal 	& 4.8 & LU & Luxembourg & 0.7  \\	
	NO &	 Norway 	& 13.7	& CH & Switzerland & 4.8 & & & \\				
	NL &	 Netherlands & 11.5	 & HU & Hungary & 4.4 & & & \\			
	BE &	  Belgium & 9.5	& DK & Denmark & 3.9	 & EU & Europe & 341.6\\	
        \bottomrule
    \end{tabularx}
    \caption{Country specific average hourly load, de-trended to their values for the year 2007.}
\label{table:one}
\end{table}

The distribution of the mismatch is shown in Figure 1 for three different countries. In all cases, $\gamma_n=1$, with three different choices for the mixing parameter ranging from fully solar ($\alpha^W_n=0.0$) to fully wind ($\alpha^W_n=1.0$). It is evident that the mismatch distributions depend on the mixing parameter $\alpha^W$, and that no mix will completely eliminate the positive and negative mismatches. Our optimization objective is to minimize the average residual load with respect to the mixing parameter $\alpha^W_n$. Since at $\gamma_n=1$ the average residual load $\langle \Delta^-_n(t) \rangle = \langle \Delta^+_n(t) \rangle$ is identical to the average excess power, this is equivalent to the minimization of average excess generation. The optimal mixing parameters are shown in Figure 2(a). Compared to the average optimal mix $\alpha_n^W \approx 0.71$, southern countries have a slightly smaller and northern countries a slightly larger value.

\begin{figure}[]
\centering
\includegraphics[width=0.99\linewidth]{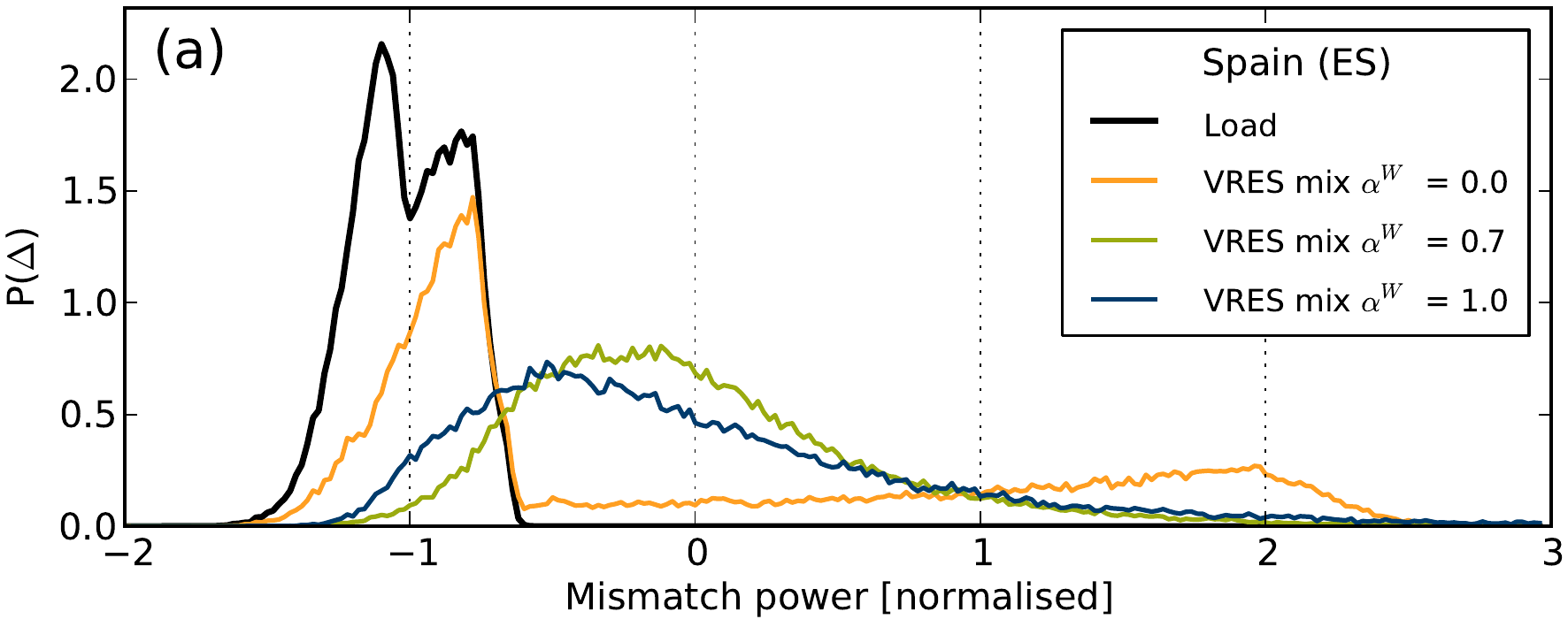}
\vspace{0.1cm}
\includegraphics[width=0.99\linewidth]{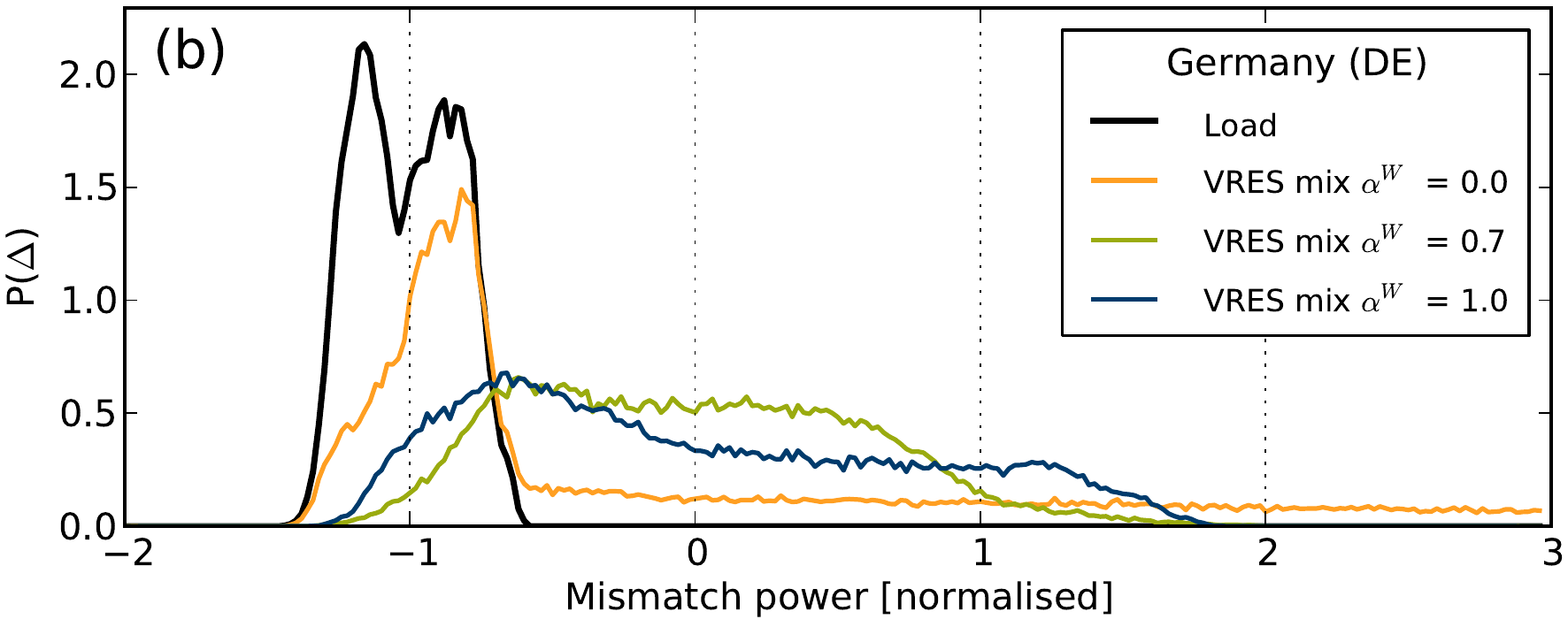}
\vspace{0.1cm}
\includegraphics[width=0.99\linewidth]{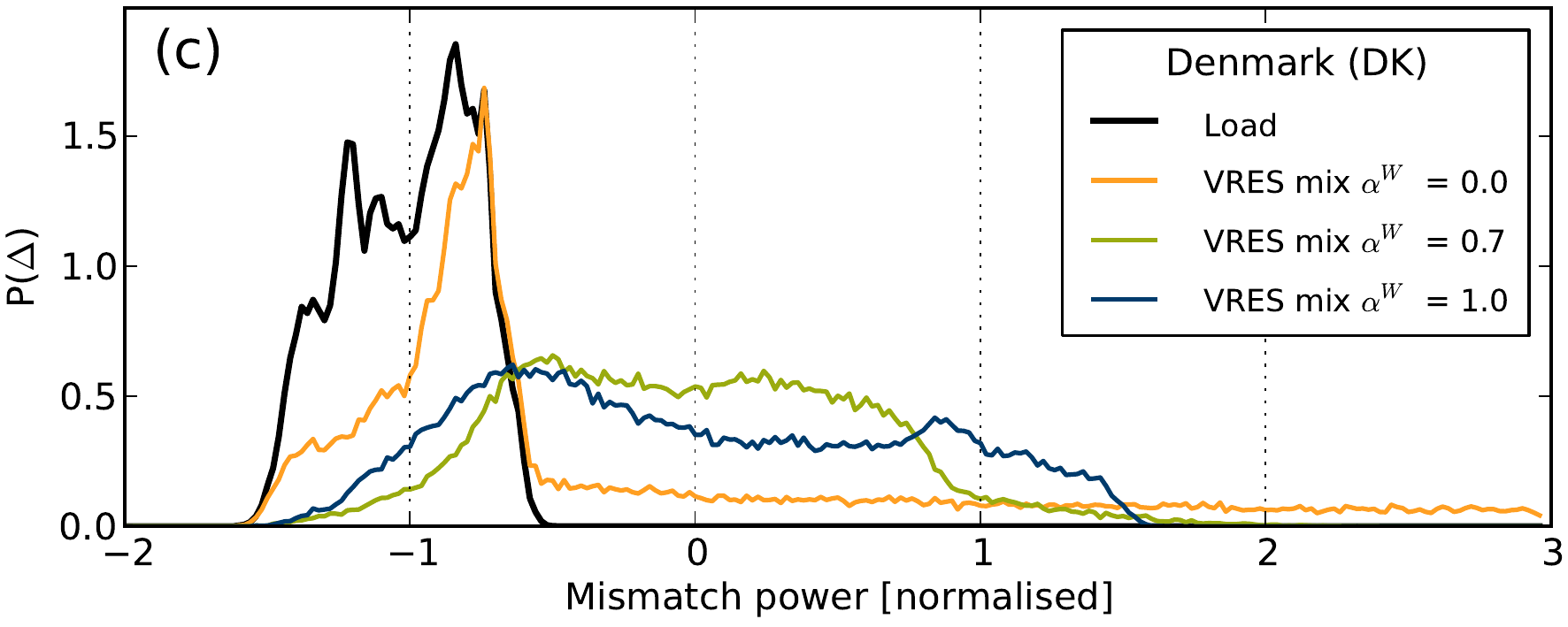}
\caption{
Mismatch distributions for (a) Spain, (b) Germany and (c) Denmark with no cross-border transmission, for $\gamma_n=1$. Values are normalized to the country specific mean load (see Table 1). The different colours represent the mixing parameters $\alpha^W=0.0, 0.7, 1.0$. For comparison, the distribution of the normalized load is also shown. 
}
\label{fig:1} 
\end{figure}

For an average country, the minimized average residual load amounts to 24\% of its average load (see Figure 2(b)). This implies that although the average wind and solar power generation are equal to the average load, 24\%  of it is generated at the wrong time. This amount has to be covered by dispatchable sources.

\begin{figure}[]
\centering
\includegraphics[width=0.99\linewidth]{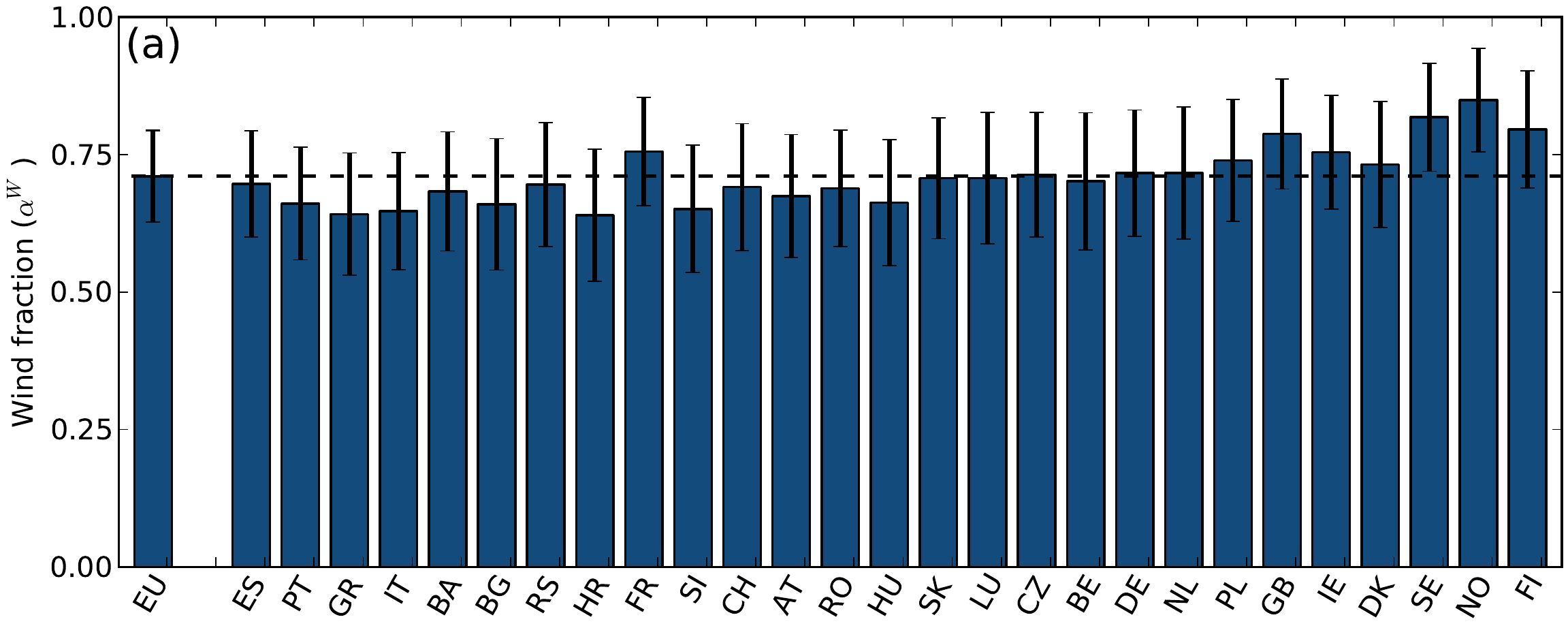}
\includegraphics[width=0.99\linewidth]{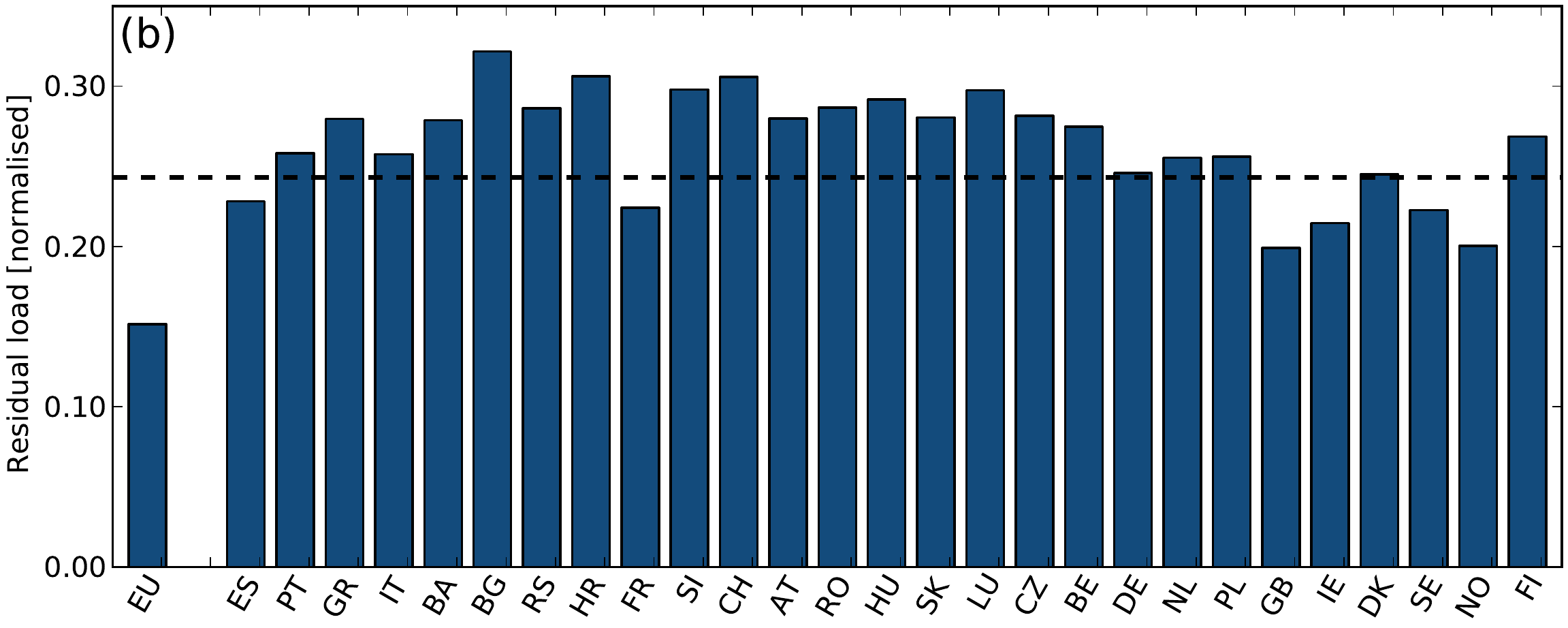}
\includegraphics[width=0.99\linewidth]{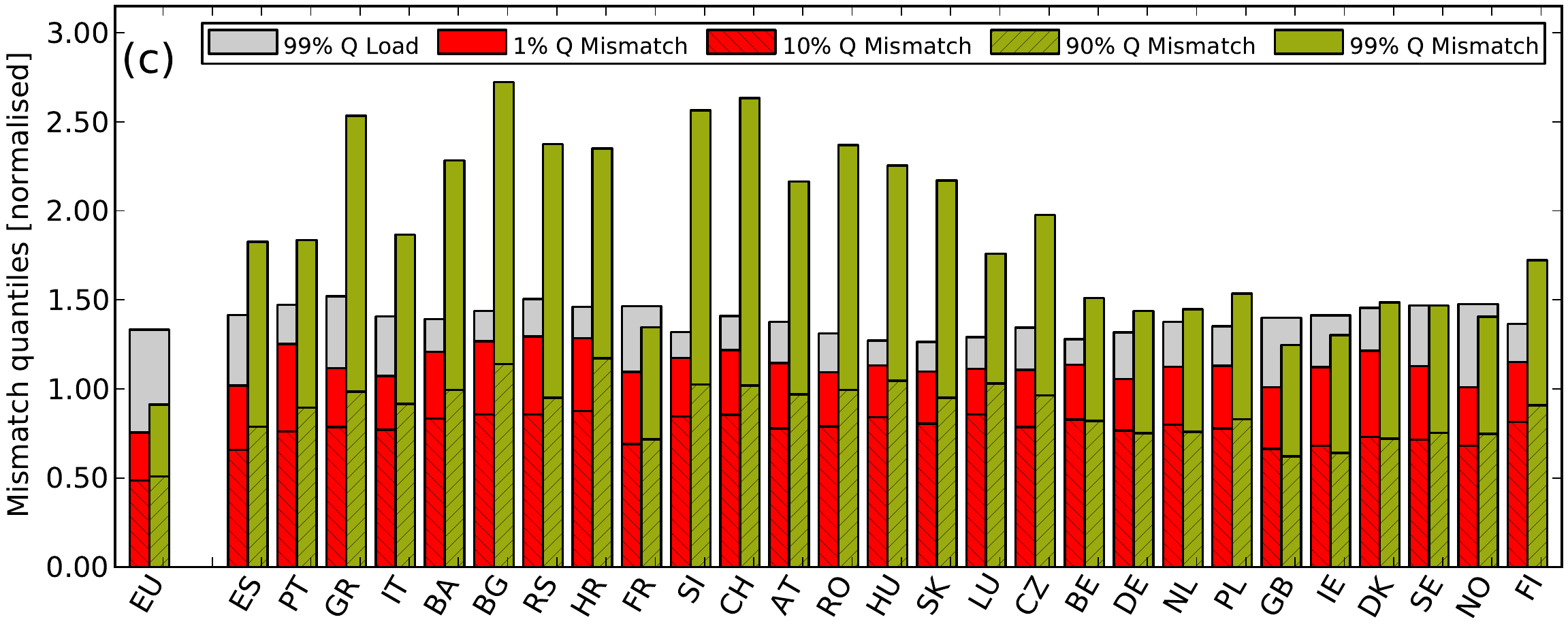}
\caption{
(a) Optimal mix minimizing the average residual load calculated independently for each country. Countries are ordered by the latitude of their geometric centre. Error bars show mixes that produce average residual loads which are larger by 1\% compared to the optimum. The dashed line indicates the optimal mix averaged over all countries and corresponds to the mix for an all-European aggregation (``EU'').
(b) Minimum average residual load for each country, obtained with the optimal mixes from (a), in units of the country specific average hourly load (see Table 1). The dashed line indicates the residual load averaged over all countries. 
(c) The 1\% (red), 10\% (striped red), 90\% (striped green), 99\% (green) quantiles of the mismatch time series based on the optimal mix from (a). As a reference the 99\% quantiles (grey) of the load time series are also indicated.}
\label{fig:2} 
\end{figure}

The 1\% and 10\% quantiles of the mismatch are shown in Figure 2(c). They give an indication of the required balancing power capacities, resulting from the residual load time series defined by (2). The 90\% and 99\% quantiles are also shown, and can be interpreted as the curtailment power, resulting from the excess power time series defined by (3). Whereas the 90\% quantiles are only a little larger than the 10\% quantiles, the 99\% quantiles turn out to be significantly larger than the 1\% quantiles. For the larger fraction of the countries the 99\% mismatch quantile is also larger than the 99\% load quantile. These are the expected results when looking again at the asymmetry of the mismatch distributions of Figure 1. 

We now compare the results for individual countries to those of the aggregated EU. The latter is assumed to have an unconstrained transmission between all countries, so that all excess generation can be shared with other countries. When each country has its optimal mix of renewables installed, the average residual load for all individual countries amounts to 24\% of the average load, shown as the dotted line in Figure 2(b). However, once all countries are aggregated, each with its own optimal mix, the resulting average residual load turns out to be $\langle \Delta_{\rm EU}^- \rangle = 15\%$. This is in full agreement with the results found in \cite{Heide:2010fk}. This means that the largest reduction in its need for balancing energy that the average country can expect from its embedding into Europe is of the order of 38\%. This sets an upper bound on what an ideal transmission system can do to reduce the need for balancing energy in Europe.

\subsection{Unconstrained power flow}

We now determine how much transmission capacity is needed for a layout to behave as though it was unconstrained. We apply the problem, as defined in equation (\ref{eq:uncprob}), to a network representing 27 European countries and the links between them. The topology of the network is based on the layout reported by ENTSO-E for winter 2010-2011 \cite{ENTSOEcaps}, and is initially assumed to have no capacity constraints; see Figure 3. Additional links in existence at the time, but with no reported NTC in \cite{ENTSOEcaps}, such as the baltic cable between Germany and Sweden and the NorNed link between Norway and the Netherlands, were not included. 

\begin{figure}[]
\centering
\includegraphics[width=0.45\linewidth]{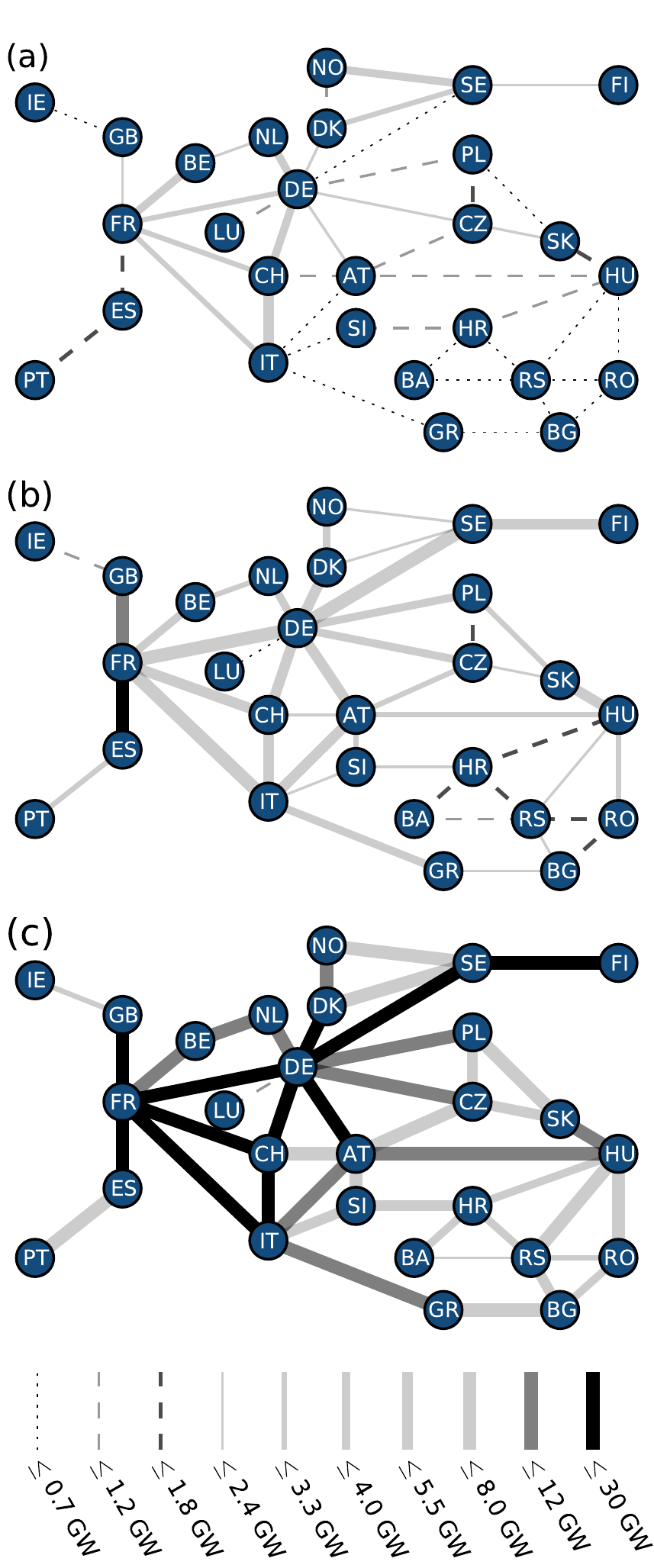}
\caption{
Transmission network topology and link capacity, with the links reported by ENTSO-E for winter 2010-2011\cite{ENTSOEcaps}. (a) Present layout capacities from winter 2010-2011. (b) Intermediate layout, with a total capacity 2.3 times larger. (c) 99\% Quantile layout, with 5.7 times the total capacity of (a). All three layouts are described in detail in Table 3. Line thickness represents the larger NTC of the interconnector.
}
\label{fig:3} 
\end{figure}

Power flows and local balancing were calculated for every hour in the eight year span, assuming the country-specific optimal mixes $\alpha^W_n$. The distribution of the resulting non-constrained power flow along a selected link is shown in Figure 4. The maximum unconstrained power flow from France to Spain amounts to 38 GW, which is larger than the combined average loads of Spain and Portugal by a factor 1.3. In the other direction the maximum flow is 75 GW, which is 1.5 times the mean load in France. 

\begin{figure}[]
\centering	
\includegraphics[width=.99\linewidth]{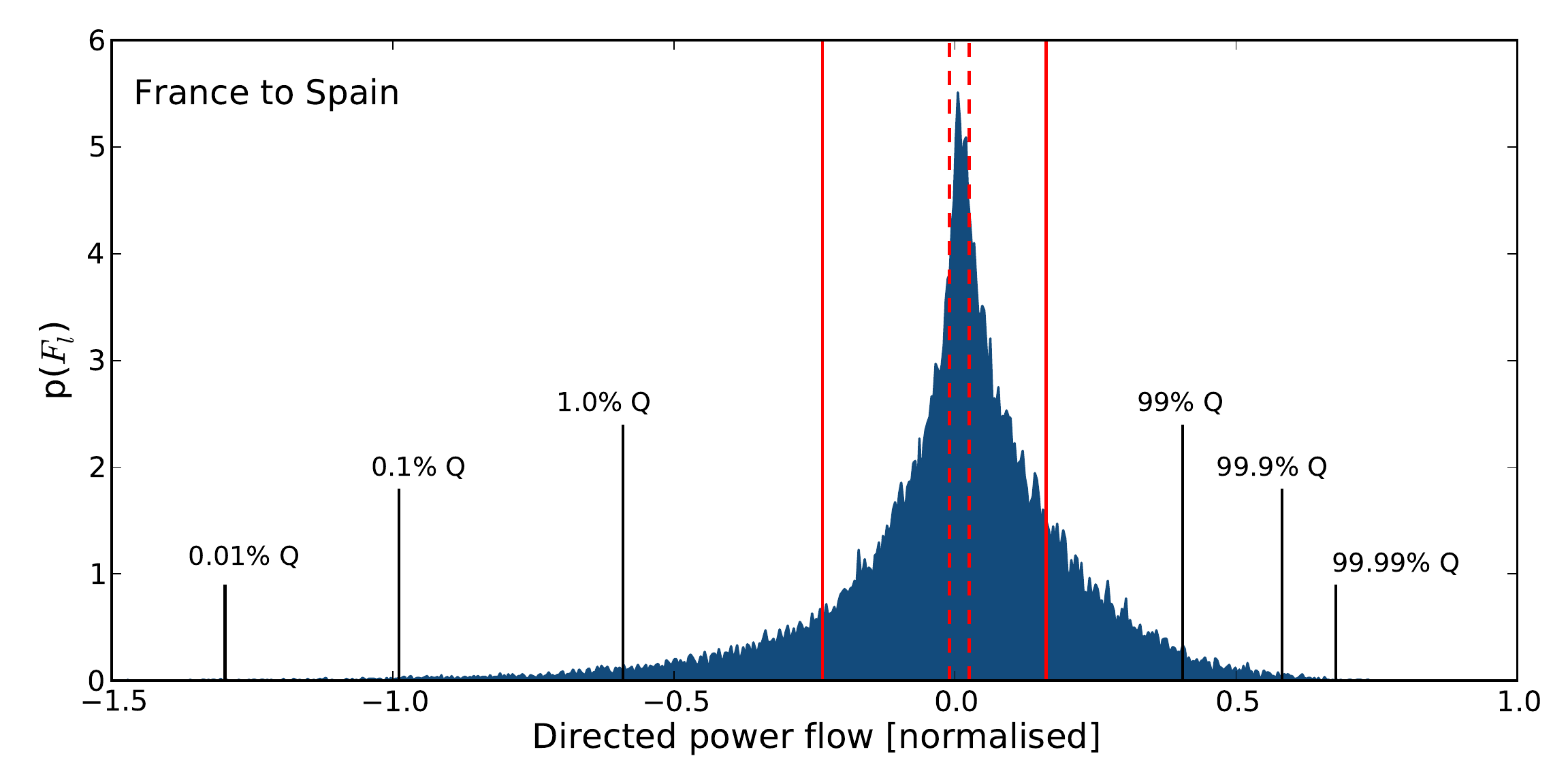}
\caption{
Distribution of unconstrained, non-zero power flows between France and Spain, normalized to the mean load in France (see Table 1). Several low- and high-quantiles are marked for illustration. The dashed red lines represent current capacities. The solid red lines show capacities as defined by the Intermediate layout. Zero-flow events occur around 46\% of the time, and are not shown.}
\label{fig:4} 
\end{figure}

We define the Unconstrained layout as that in which all links have capacities equal to the maximum recorded exchange, so that power can flow unconstrained along the interconnectors. By construction, these give rise to the full benefit of cooperation, as they allow the interactions between countries to be identical to one in which they are all aggregated. The sum 
\begin{equation}
T=\sum_{l=1}^L \max \left \{ |f_l^-| , |f_l^+| \right \}
\end{equation}
over the larger NTC value of each interconnector in the present layout, as reported by ENTSO-E for winter 2010-2011\cite{ENTSOEcaps}, adds up to around 69 GW. With 831 GW the Unconstrained layout capacities are twelve times larger. These unconstrained capacities are determined by single, one-hour events over eight years of data. Therefore, we consider the 1\% and 99\% quantiles of the flow distributions to define a reduced, directed capacity layout, which we call the 99\% Quantile layout; see again Figure 4. This means that power will flow unobstructed for 98\% of the time. The remaining 2\% corresponds to around one week per year. The 99\% Quantile layout comes with 392 GW in total and is roughly half as large as the Unconstrained layout, but still 5.7 times larger than today's interconnector capacities. See again Table 3. 

\subsection{Constrained power flow}

To determine what fraction of the benefit of transmission is obtained with a non-ideal, limited transmission capacity, we deal with constrained power flows as defined in (\ref{eq:conprob}). This allows the determination of a compromise between the reduction in balancing energy and the increase in total transmission capacity. 

\begin{table}[t]
\scriptsize
\centering
\begin{tabularx}{133mm}{l S c S S}
	\toprule
      	{\bf Layout}    & {Total transmission } & {$E_B$} & {$E_B$ (percentage of} &  {Benefit $\beta$ of } \\
	{} & {capacity [GW]} & {[TWh]} &  {annual consumption)} &{transmission}\\
    	\midrule
	\rowcolor[gray]{0.95}Zero transmission & 0 & 727 & 24.3\% & 0.0\%\\
	Present layout & 69 & 633 & 21.2\% & 34.2\%\\
	\rowcolor[gray]{0.95}Intermediate layout & 157 & 535 & 17.9\% & 69.8\%\\
	99\% Quantile layout & 392 & 461 & 15.4\% & 96.7\%\\
	\rowcolor[gray]{0.95}Unconstrained layout& 831 & 452 & 15.1\% & 100.0\%\\
	\bottomrule	
\end{tabularx}
\caption{Studied layouts with their total installed transmission capacities, the required total European balancing energy $E_B$ in absolute and relative terms, and the benefit of transmission $\beta$.}
\label{table:two}
\end{table}

As can be seen in Table 2, the 99\% quantile capacities provide, with $\beta = 96.7\%$, most of the benefit of the Unconstrained layout with less than half of the total installed capacity. The layout defined by these 99\% quantiles can be seen in Figure 3(c). It is also noteworthy that today's capacities already provide 34.2\% of the benefit of transmission, if applied to this scenario. In order to find out how the benefit scales with increasing transmission capacities, ways of interpolating between today's system and the larger layouts must first be defined.

Interpolation A is an upscaling of present capacities with a linear factor $a$. That is, for a directed link $l$, the limits are defined by 
\begin{equation}
f_l^{\rm A} = \min \left \{ a f_l^{\rm today} , f_l^{99 {\rm \% Q}} \right \}
\; ,
\end{equation}
where $f_l^{\rm today}$ represents the NTC of the link as of winter 2010-2011 and $f_l^{99 {\rm \% Q}}$ those of the 99\% Quantile layout.

Interpolation B involves a linear reduction of the 99\% quantile capacities with factor $b$, that is
\begin{equation}
f_l^{\rm B} = b f_l^{99 {\rm \% Q}}
\; .
\end{equation}
Interpolation C defines the capacity layout
\begin{equation}
f_l^{\rm C} = f_l^{c {\rm Q}}
\; ,
\end{equation}
which allows unconstrained flow for a percentage $c$ of time, as shown by the different quantiles in Figure 4. Here, more capacity is allocated to more transited links than to less used ones.

Figure 5 shows the average balancing energy required by all nodes for different transmission layouts, following all three interpolations. We use the larger of the NTC on each link direction as a proxy to estimate the total installed capacity shown on the x-axis.

\begin{figure}
\centering
\includegraphics[width=0.6\linewidth]{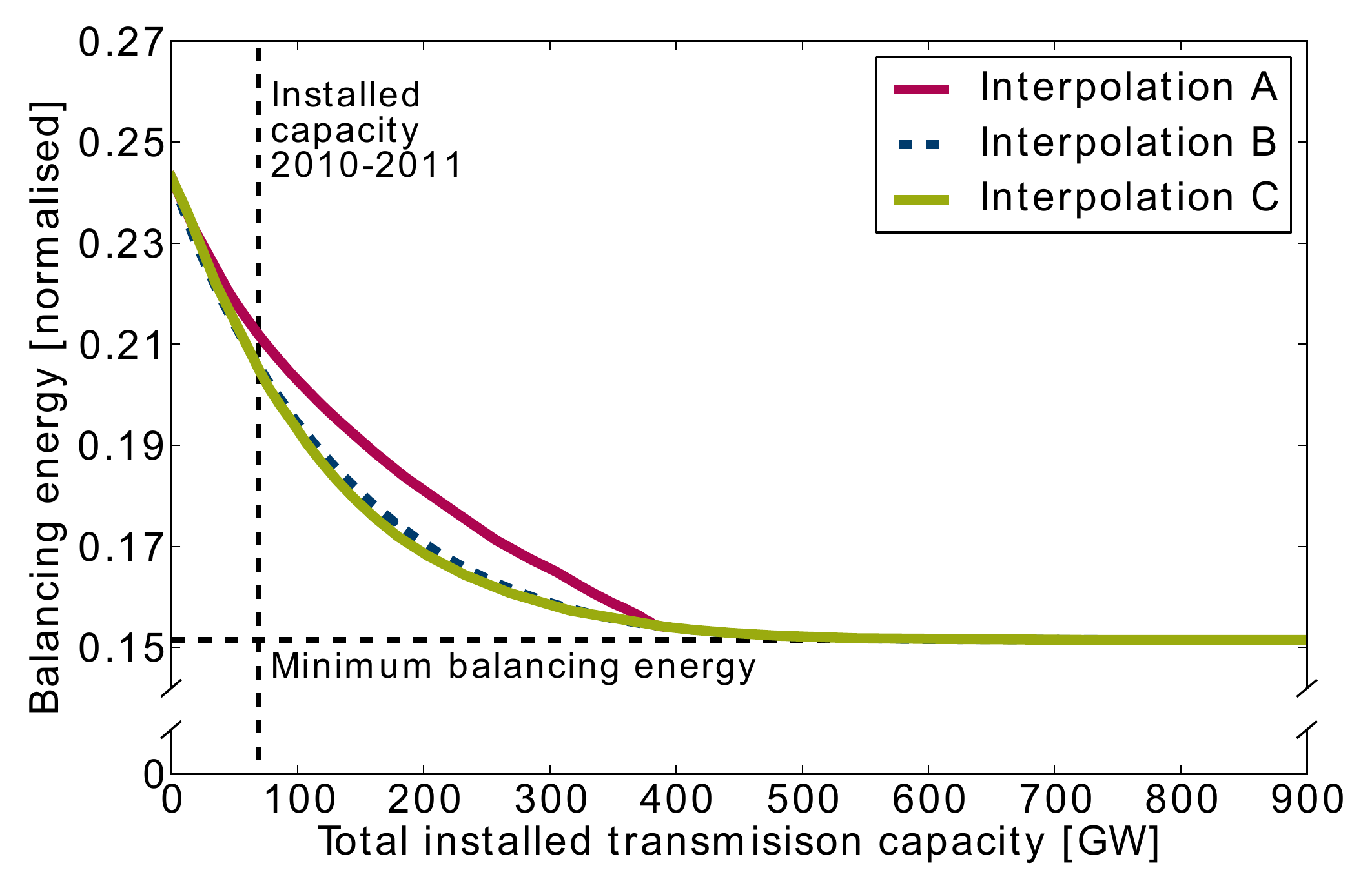}
\caption{
Balancing energy as a function of the total installed transmission capacity for different interpolations, normalized to the total annual consumption. The vertical dashed line indicates the total installed transmission capacity as of 2010-2011, and the horizontal dashed line the total balancing required for the unconstrained flow. Notice that the intersection at $T=0$ is congruent with the results shown in Figure 2 for the average balancing required. The large $T$ asymptotics for the more efficient scenarios is congruent with the balancing requirement for the aggregated European countries, also shown in Figure 2. 
}
\label{fig:5} 
\end{figure}

Today's existing transmission capacity already provides a reduction in balancing energy from 25\% to 21\% of the annual consumption, providing 34\% of the benefit of transmission. This is not as much as a more efficient distribution of the resources could have achieved. With the same total amount of transmission capacities, layouts B and C lead to a balancing energy of 20\% of the annual consumption, with the benefit of transmission being 45\%. An increase of today's total capacities by a factor of 2.3 in the way described by interpolations B or C, from 69 GW to 157 GW, will double the benefit of transmission to 70\%, reducing the balancing energy to about 18\% of the annual consumption. This new layout is between today's transmission capacities and the capacities defined by the 99\% quantiles. Hence, we denote it as the Intermediate layout. It can be seen in Figure 3(b), and its capacities are also listed in Table 3. Overall results for the Intermediate layout are presented in Table 2.

Notice that some of the capacities under the Intermediate layout are smaller than the ones in the present one. The link between Norway and Sweden, for instance, is $f_{\rm NO-SE}^{\rm today}=3.90$  GW while $f_{\rm NO-SE}^{\rm inter}=1.99$   GW. This is due to the fact that the Intermediate layout is presented as a reduction of the 99\% quantiles, and not as an expansion of the present one.

\section{Discussion: country perspectives}

After having quantified how much the strengthening of the interconnectors between the countries reduces the total balancing energy of Europe, we are now interested in what the different transmission capacity layouts imply for the import / export capabilities and the balancing power capacities of the single countries.

\subsection{Limits to export and import capabilities}

Let us have a look again at Table \ref{table:two}. Compared to the Zero transmission layout, the Unconstrained layout reduces the European balancing energy by 37.8\%, from 727 TWh to 452 TWh. It is caused by the power flows from countries with an excess power generation to those with a deficit. The 37.8\% can also be interpreted as the maximum capability for imports and exports. This number can not become larger since it is already based on the unconstrained transmission layout. It draws the limit to the benefit that geographical dispersion of VRES can bring to Europe. For the constrained transmission layouts the import / export capabilities are smaller. The relative reduction in balancing energy is 12.8\%, 26.3\% and 36.6\% for the present, the intermediate and the 99\% quantile transmission layout, respectively.  

\begin{figure}[]
\centering
\includegraphics[width=0.99\linewidth]{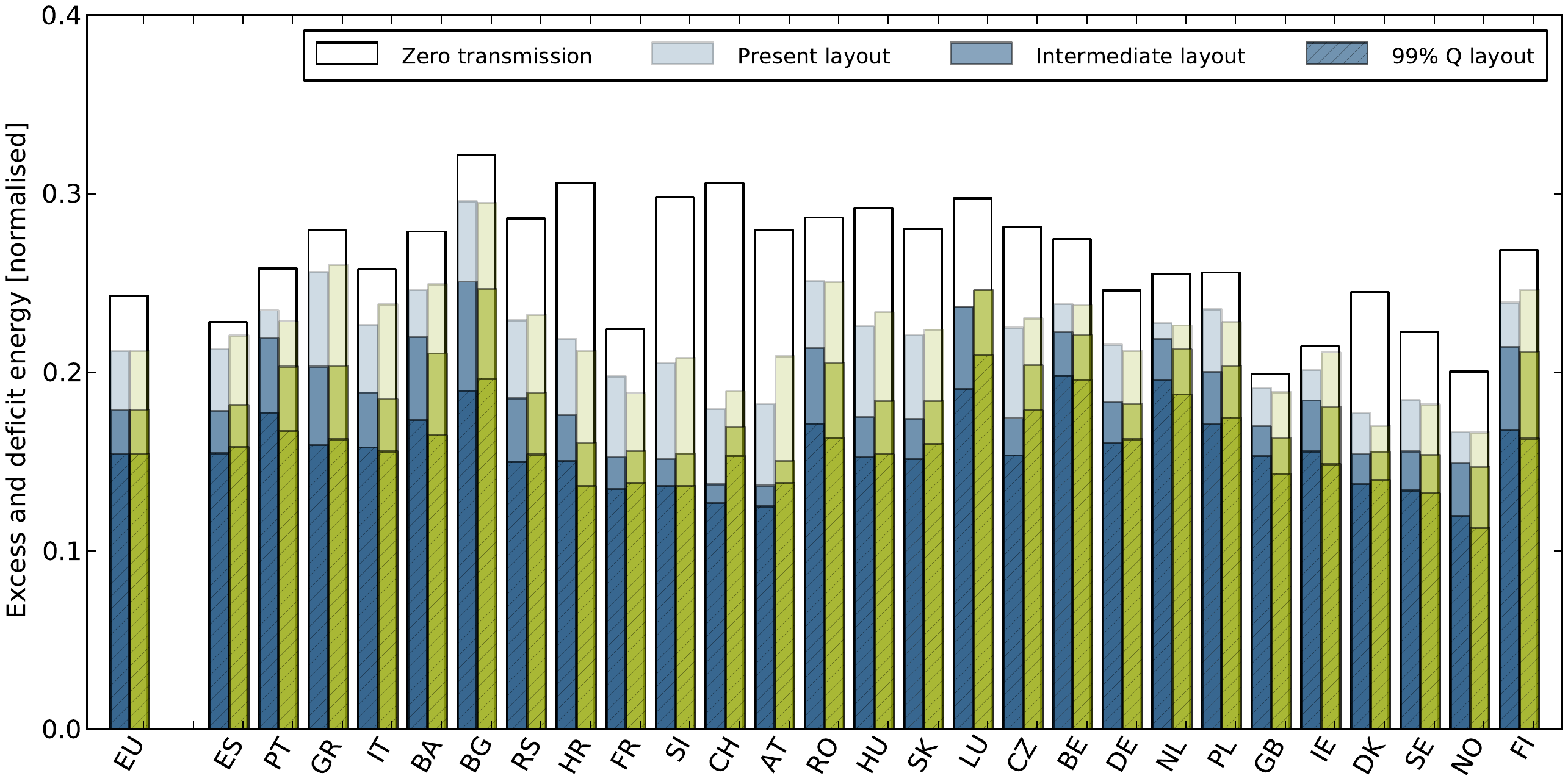}
\caption{
Dependence of the country-specific average residual loads (blue) and excess powers (green) on the Zero transmission layout, the Present layout, the Intermediate layout and the 99\% Quantile layout. Residual loads and excess powers are normalized by the average load. Bars are not stacked. The optimal mixes $a^W_n$ from Figure 2 (a) have been used for each country. Note, that the residual loads for the Zero transmission layout are identical to those shown in Figure 2 (b). The EU bar represents four points on the "interpolation C" curve of Figure 5. 
}
\label{fig:6} 
\end{figure}

So far the reduction in balancing energy has been discussed for total Europe only. For the single countries the reduction does not need to be the same, although the average over all countries has to reproduce Europe's reduction in balancing energy. Figure \ref{fig:6} illustrates the time-averaged country-specific residual loads for the different transmission capacity layouts. Apparently some countries have better import capabilities than others. Most of the small middle and southern countries show a pronounced reduction of balancing energy from the present to the intermediate transmission capacity layout. As to the larger countries, France benefits the most from imports, and Great Britain the least. Compared to the present transmission capacity layout, the intermediate transmission layout reduces the balancing energies between 4\% (for the Netherlands) and 26\% (for Slovenia). For the 99\% quantile transmission layout the respective reductions are in the range between 14\% (for the Netherlands) and 38\% (for Greece).

Figure \ref{fig:6} also shows the country-specific average excess powers for the various transmission capacity layouts. For some countries like Great Britain, Ireland, the Netherlands, Romania, Croatia and Portugal the average excess powers turn out to be smaller than their average residual loads. With the further extension of the present transmission capacity layout, these countries will be able to export more than what they will be able to import. Other countries like Switzerland, Austria and Czech Republic will become strong importers. It is interesting to note that the stronger transmission capacity layouts quite naturally lead to stronger import / export imbalances for the countries.

\subsection{Balancing power capacities}

\begin{figure}
\centering
\includegraphics[width=0.99\linewidth]{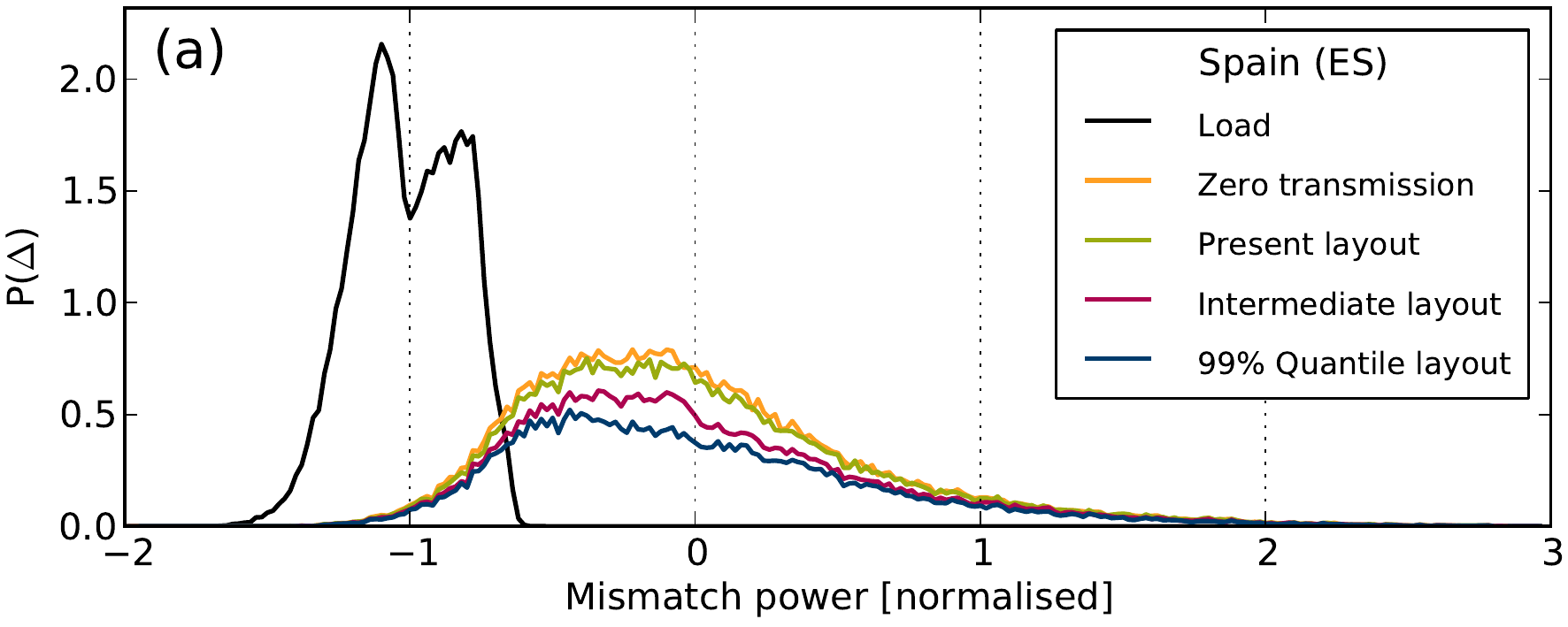}
\vspace{0.1 cm}
\includegraphics[width=0.99\linewidth]{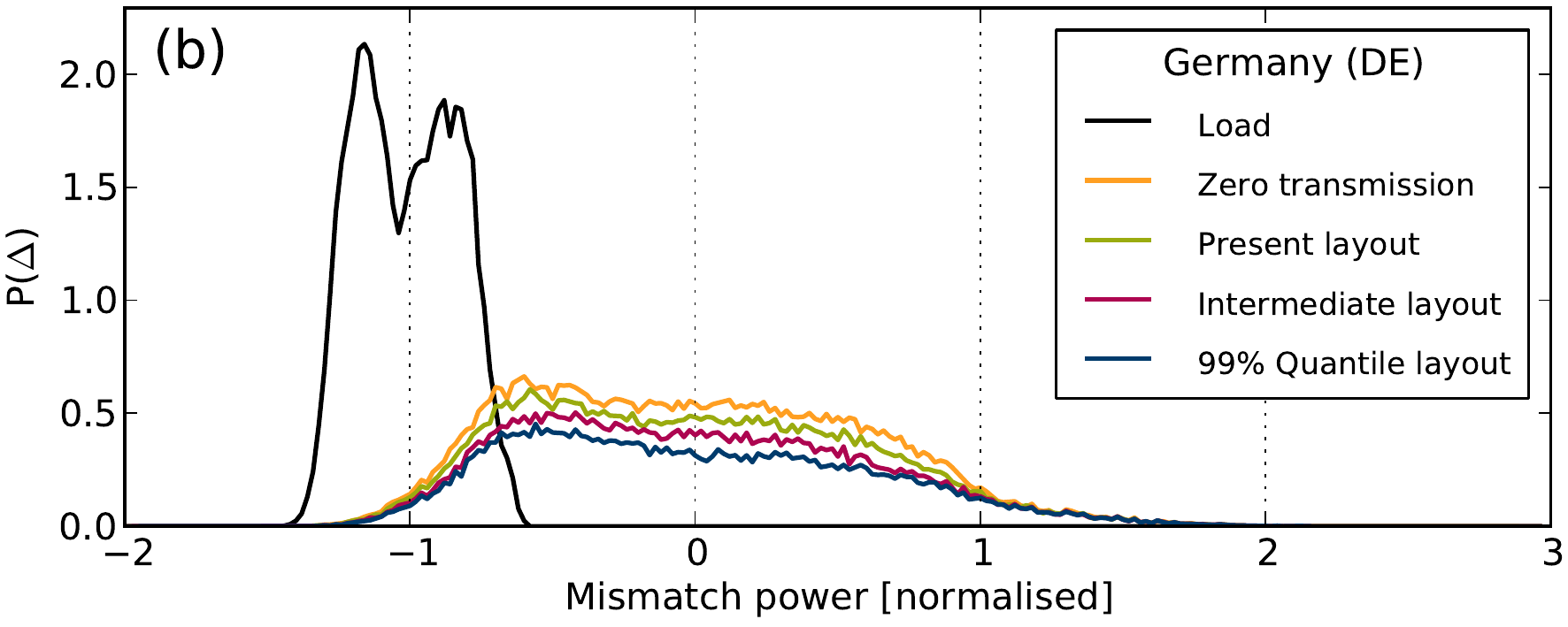}
\vspace{0.1 cm}
\includegraphics[width=0.99\linewidth]{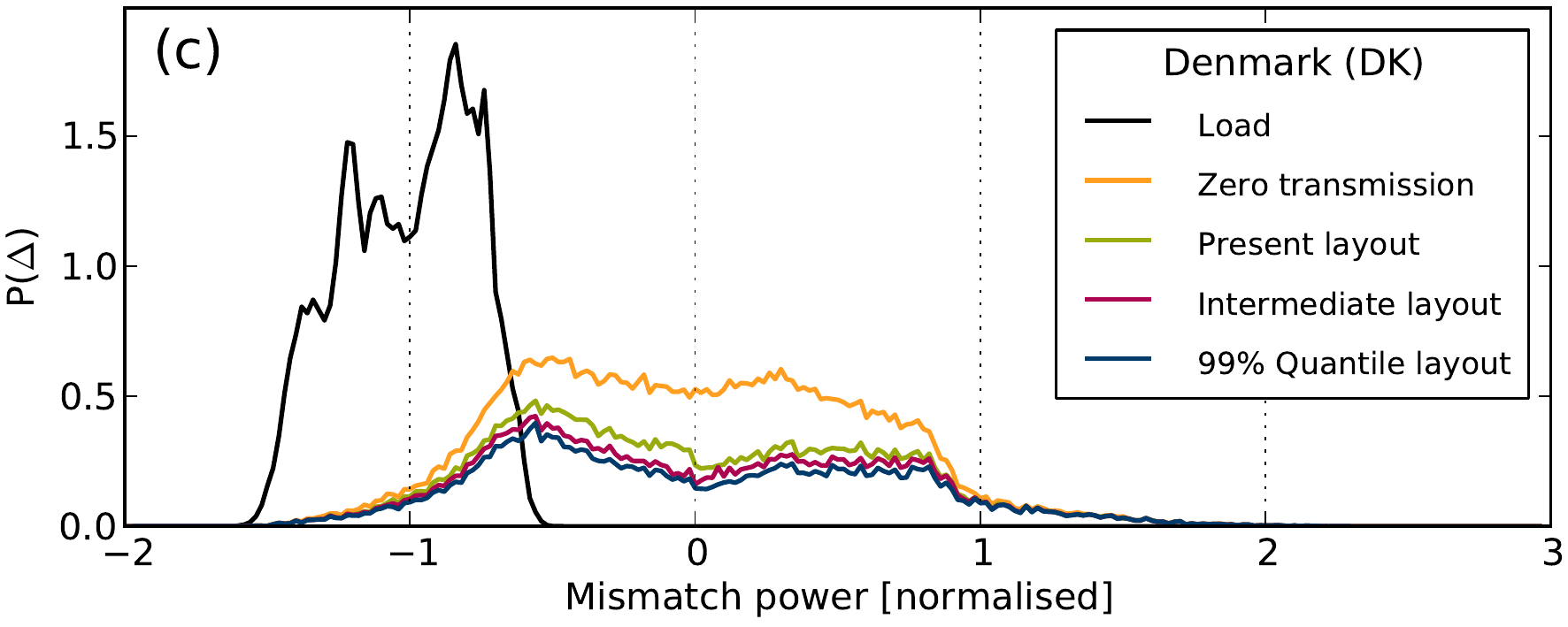}
\caption{
Normalized non-zero mismatch distributions for (a) Spain, (b) Germany and (c) Denmark, with a country-specific optimal mix of wind and solar. Apart from the zero transmission scenario (orange), three layouts are shown: the Present layout (green), the Intermediate layout (purple) and the 99\% Quantile layout (blue). For comparison, the distribution of the normalized load (black) is also shown.
}
\label{fig:7} 
\end{figure}

Figure \ref{fig:7} shows the distributions of non-zero mismatches for three selected countries, as they change for different transmission capacity layouts. The central part of the distributions are lowered as the transmission capacity layouts become stronger. This connects nicely to the results on the reduction of average residual load and excess power, which have been discussed in the previous Subsection. Contrary to the central part, the tails of the distributions remain unaffected by the increase in transmission capacities. Referring to Figure \ref{fig:2}(c), this implies that the 99\% quantiles of the residual load and the excess power for the non-zero transmission capacity layouts are not reduced when compared to those for the zero transmission capacity layout. In other words, it appears that increased transmission does not affect the maximum balancing power, and that Europe has to keep its present dispatchable power generation capacity. However, this would be a hasty conclusion.

\begin{figure}
\centering
\includegraphics[width=0.99\linewidth]{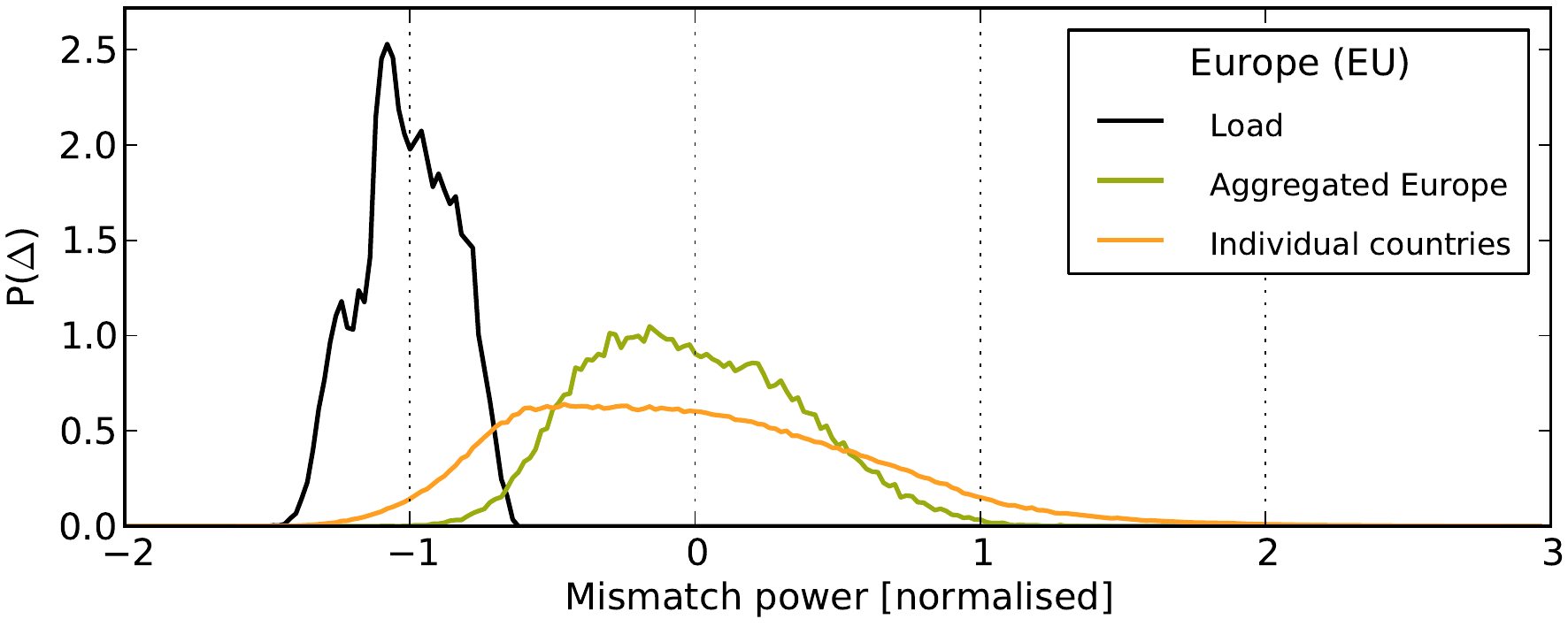}
\vspace{0.1 cm}
\caption{
Normalized non-zero mismatch distribution for Europe resulting from the unconstrained transmission scenario (green) and for an average country in the zero transmission scenario (orange). For comparison, the distribution of the normalized European load (black) is also shown.
}
\label{fig:8} 
\end{figure}

The objective of the current power flow modelling, as presented in Section 2, has been first to reduce the overall balancing energy the most, and then to determine the most localized power flow across the network. The objective has not been to reduce the high quantiles of the residual load at each node of the network. In other words, the reduction of balancing power capacities has not been prioritized in the current power flow modelling. Figure \ref{fig:2}(c) gives an indication about how much the balancing power capacities can be reduced in principle. The first bar of this figure shows the 1\% quantile of the mismatch for an aggregated Europe. With 73\% of the average European load it is significantly smaller than 108\% which holds true for an average independent single country. This is also nicely visualized in Figure \ref{fig:8}, which shows the respective mismatch distribution for Europe with unconstrained transmission and for an average country in the zero transmission scenario. Compared to the single-country distribution, both tails of the European distribution are significantly shifted towards zero mismatch. The explanation for this reduction in the overall balancing capacities lies in the sharing of balancing capacities between countries. 

This effect of shared balancing is not taken into account in the current power flow modelling, where after import of renewable power the countries balance their remaining deficit fully by themselves. Shared balancing implies that countries are allowed to balance part of their remaining deficit from balancing capacities of other countries, which are not fully used at the same time. This of course also leads to additional power flows. A self-consistent treatment of power-flow modelling with shared balancing capacities will not be given here. Other ways of reducing peak balancing power could include demand flexibility, the use of energy storage, and penetration of renewables larger than 100\%. For the moment we will leave all of this open for future investigations.

\section{Conclusion}
\label{sec:five}

We have quantified the benefit of the far-future pan-European transmission system, when by assumption all countries have reached a 100\% penetration of combined wind and solar power generation. Two extreme transmission scenarios can be treated without any power flow calculations. For the scenario with zero inter-connectors, the countries can be discussed separately. Their annual balancing energy, which is required to cover the negative mismatches between their renewable power generation and their load, depends on the mixing parameter $\alpha^W$ between wind and solar power generation, becomes minimal at $\alpha^W \approx 0.71$ (with a small dependence on the geographical latitude), and then results to be around 24\% of their annual electricity consumption. For the scenario with infinitely strong inter-connectors and no transmission losses, all country-specific mismatches can be directly aggregated into an overall European mismatch. This leads to a required balancing energy, which is only 15\% of the total annual European electricity consumption. The difference 24\% - 15\% between the outcomes of the two extreme scenarios represents an upper limit on the benefit that a pan-European transmission system can provide. In other words, no transmission layout is able to further reduce the required balancing energy. 

In order to estimate the benefit of transmission capacities constrained between zero and infinity, a novel modelling approach has been presented, which calculates the constrained power flows for a pan-European network. Based on highly resolved spatio-temporal weather data, it first minimizes the all-European residual load, so that as much of the overall renewable power generation is used as possible. In a second step, the square sum over all flows is minimized, leading to a very local flow pattern and avoiding additional assumptions regarding market-economic policies. This novel modelling approach reveals that an infinitely strong European transmission network suffices to be twelve times as strong as today's total inter-connector capacities. A capacity layout six times as large as today's provides 97\% of the benefit of transmission, which is almost as good as the infinitely strong layout for the reduction of required balancing energy. For weaker transmission networks, the relationship between the need for balancing energy and the total capacity of the transmission layout turns out to be non-linear and convex. A good compromise between these two conflicting objectives, i.e.\ on the one hand reducing the balancing energy as much as possible and on the other to increase the transmission capacities as little as possible, appears to be an Intermediate layout with total transmission capacities being 2.3 larger than today's capacities. This Intermediate layout leads to a 70\% benefit of transmission. Compared to the Zero transmission layout it reduces the required European balancing energy from 24\% to 18\% of the annual consumption. This reduction also implies that an average country participating in the pan-European transmission network can only import 26\% of its own balancing needs, and that the remaining 74\% have to come from own balancing resources.

The presented findings have focused on the transmission needs in a fully renewable European power system. The penetration of combined wind and solar power generation has been assumed to be exactly 100\%. Of course, this can be generalized to penetrations below 100\%, and thus provide information on the required ramp-up of transmission needs all the way from today's penetration to the far-future 100\% penetration \cite{Becker2012}. In this respect, also a modification of the network topology should be discussed, like for example the addition of a few new selected links across the North Sea and of a few more countries. In order to reduce the need for balancing energy further, maybe even below the discussed 15\% limit of the annual 
electricity consumption, it will be interesting to look at renewable penetrations above 100\% and discuss more synergies between balancing, transmission and storage.

\bigskip
\noindent {\bf Acknowledgements:} The authors thank Uffe Poulsen and Morten G.\ Rasmussen for insightful discussions. S.B. gratefully acknowledges financial support from O. and H. St\"ocker, and G.B.A. financial support from DONG Energy and The Danish National Advanced Technology Foundation. 
\section*{Bibliography}
\label{sec:bib}

\bibliographystyle{unsrt}
\bibliography{references}
\newpage

\section{Appendix}
\begin{table}[h!]
	\scriptsize
	\centering
    \begin{tabularx}{{132mm}}{  l c c c c }
    \toprule
    \hiderowcolors
    Link & Present Layout & Intermediate Layout & 99 \% Q Layout & Unconstrained Layout \\ 
     & (GW) & (GW) & (GW) & (GW)\\
    \hline
    \showrowcolors
    
 &  0.47  &  2.39  &  5.97  &  15.32 \\
\multirow{-2}{*}{ AT  $\rightleftarrows$  CH } &  1.20  &  2.21  &  5.53  &  12.87 \\
\rowcolor{rgrey} &  0.60  &  3.09  &  7.72  &  14.56 \\
\rowcolor{rgrey} \multirow{-2}{*}{ AT  $\rightleftarrows$  CZ } &  1.00  &  2.76  &  6.9  &  12.33 \\
 &  2.00  &  5.42  &  13.56  &  30.81 \\
 \multirow{-2}{*}{ AT  $\rightleftarrows$  DE } &  2.20 &  4.50  &  11.25  &  18.15 \\
\rowcolor{rgrey} &  0.80  &  2.37  &  5.91  &  11.37 \\
\rowcolor{rgrey}\multirow{-2}{*}{ AT  $\rightleftarrows$  HU } &  0.80  &  3.26  &  8.15  &  20.36 \\
 &  0.22  &  3.44  &  8.59  &  16.52 \\
\multirow{-2}{*}{ AT  $\rightleftarrows$  IT } &  0.28  &  4.17  &  10.41  &  21.69 \\
\rowcolor{rgrey} &  0.90  &  2.00  &  4.99  &  7.97 \\
\rowcolor{rgrey} \multirow{-2}{*}{ AT  $\rightleftarrows$  SI } &  0.90  &  2.69  &  6.73  &  15.66 \\
 &  0.60  &  1.42  &  3.55  &  7.24 \\
 \multirow{-2}{*}{ BA  $\rightleftarrows$  HR } &  0.60  &  0.91  &  2.27  &  3.85 \\
\rowcolor{rgrey} &  0.35  &  0.95  &  2.37  &  4.16 \\
\rowcolor{rgrey} \multirow{-2}{*}{ BA  $\rightleftarrows$  RS } &  0.45  &  0.71  &  1.78  &  4.04 \\
&  2.30  &  3.41  &  8.51  &  16.06 \\
\multirow{-2}{*}{ BE  $\rightleftarrows$  FR } &  3.40  &  3.79  &  9.48  &  18.96 \\
\rowcolor{rgrey}  &  2.40  &  3.20  &  8.01  &  16.35\\
\rowcolor{rgrey} \multirow{-2}{*}{ BE  $\rightleftarrows$  NL } &  2.40  &  3.09  &  7.72  &  15.14 \\
&  0.55  &  2.22  &  5.55  &  15.01 \\
 \multirow{-2}{*}{ BG  $\rightleftarrows$  GR } &  0.50  &  1.97  &  4.92  &  10.93 \\
\rowcolor{rgrey} &  0.60  &  1.33  &  3.31  &  7.05 \\
\rowcolor{rgrey} \multirow{-2}{*}{ BG  $\rightleftarrows$  RO } &  0.60  &  1.04  &  2.59  &  5.28 \\
&  0.45  &  2.02  &  5.05  &  10.97 \\
\multirow{-2}{*}{ BG  $\rightleftarrows$  RS } &  0.30  &  1.28  &  3.19  &  5.85 \\
\rowcolor{rgrey} &  3.50  &  4.95  &  12.37  &  23.79 \\
\rowcolor{rgrey} \multirow{-2}{*}{ CH  $\rightleftarrows$  DE } &  1.50  &  4.40  &  11.00  &  20.01 \\
 &  1.10  &  4.45  &  11.14  &  28.27\\
 \multirow{-2}{*}{ CH  $\rightleftarrows$  FR }  &  3.20  &  4.87  &  12.19  &  28.34 \\
 \rowcolor{rgrey} &  4.17  &  3.65  &  9.12  &  15.20 \\
 \rowcolor{rgrey} \multirow{-2}{*}{ CH  $\rightleftarrows$  IT } &  1.81  &  4.86  &  12.14  &  26.99 \\
&  2.30  &  3.36  &  8.40  &  19.54 \\
\multirow{-2}{*}{ CZ  $\rightleftarrows$  DE } &  0.80  &  3.01  &  7.52  &  14.45 \\
\rowcolor{rgrey}   &  0.80  &  1.68  &  4.21  &  8.26 \\
\rowcolor{rgrey} \multirow{-2}{*}{ CZ  $\rightleftarrows$  PL } &  1.80  &  1.77  &  4.43  &  7.92\\
 &  2.20  &  1.38  &  3.45  &  5.95 \\
\multirow{-2}{*}{ CZ  $\rightleftarrows$  SK } &  1.20  &  1.81  &  4.53  &  9.54 \\
\rowcolor{rgrey} &  1.55  &  4.93  &  12.32  &  21.40 \\
\rowcolor{rgrey} \multirow{-2}{*}{ DE  $\rightleftarrows$  DK } &  2.08  &  5.37  &  13.43  &  25.83 \\

\bottomrule
    \end{tabularx}
    \caption*{Table 3.a: Interconnector capacities  for different layouts.}
\end{table}

\begin{table*}
	\centering
    	\scriptsize
    \begin{tabularx}{{132mm}}{  l c c c c }
    \toprule
    \hiderowcolors
    Link & Present Layout & Intermediate Layout & 99\% Q Layout & Unconstrained Layout \\ 
     & (GW) & (GW) & (GW) & (GW)\\
    \hline
    \showrowcolors

&  3.20  &  7.14  &  17.86  &  37.21 \\
\multirow{-2}{*}{ DE  $\rightleftarrows$  FR} &  2.70  &  7.76  &  19.39  &  42.16 \\
\rowcolor{rgrey}  &  0.98  &  0.25  &  0.63  &  0.90 \\
\rowcolor{rgrey}  \multirow{-2}{*}{ DE  $\rightleftarrows$  LU } &  NRL  &  0.35  &  0.88  &  1.51 \\
&  3.85  &  3.58  &  8.95  &  18.45 \\
\multirow{-2}{*}{ DE  $\rightleftarrows$  NL }  &  3.00  &  3.65  &  9.14  &  20.56 \\
\rowcolor{rgrey} &  1.20  &  3.18  &  7.95  &  17.27\\
\rowcolor{rgrey} \multirow{-2}{*}{ DE  $\rightleftarrows$  PL }&  1.10  &  3.31  &  8.28  &  20.21  \\
 &  0.60  &  6.69  &  16.72  &  28.34 \\
 \multirow{-2}{*}{ DE  $\rightleftarrows$  SE } &  0.61  &  7.49  &  18.72  &  36.00 \\
\rowcolor{rgrey}  &  0.95  &  2.82  &  7.05  &  11.55 \\
\rowcolor{rgrey} \multirow{-2}{*}{ DK  $\rightleftarrows$  NO }  &  0.95  &  3.32  &  8.29  &  15.49\\
 &  2.44  &  1.94  &  4.84  &  7.44  \\
\multirow{-2}{*}{ DK  $\rightleftarrows$  SE }  & 1.98  &  2.32  &  5.80  &  10.07\\
\rowcolor{rgrey}   &  0.50  &  12.07  &  30.18  &  75.44\\
\rowcolor{rgrey} \multirow{-2}{*}{ ES  $\rightleftarrows$  FR } &  1.30  &  8.29  &  20.72  &  37.29 \\
 &  1.70  &  1.87  &  4.67  &  7.79 \\
\multirow{-2}{*}{ ES  $\rightleftarrows$  PT } &  1.50  &  2.42  &  6.05  &  12.35 \\
\rowcolor{rgrey} &  1.65  &  4.90  &  12.25  &  20.76 \\
\rowcolor{rgrey} \multirow{-2}{*}{ FI  $\rightleftarrows$  SE } &  2.05  &  3.48  &  8.70  &  11.92 \\
 &  2.00  &  9.20  &  23.00  &  47.44 \\
\multirow{-2}{*}{ FR  $\rightleftarrows$  GB } &  2.00  &  10.43  &  26.09  &  42.34 \\
\rowcolor{rgrey} &  2.58  &  6.87  &  17.18  &  39.95 \\
\rowcolor{rgrey} \multirow{-2}{*}{ FR  $\rightleftarrows$  IT } &  0.99  &  7.63  &  19.08  &  42.93 \\
 &  0.45  &  1.03  &  2.56  &  4.45 \\
\multirow{-2}{*}{ GB  $\rightleftarrows$  IE } &  0.08  &  0.90  &  2.26  &  4.34 \\
\rowcolor{rgrey} &  0.50  &  3.98  &  9.96  &  20.32 \\
\rowcolor{rgrey} \multirow{-2}{*}{ GR  $\rightleftarrows$  IT } &  0.50  &  2.72  &  6.79  &  11.52 \\
 &  0.80  &  1.45  &  3.62  &  7.10 \\
\multirow{-2}{*}{ HR  $\rightleftarrows$  HU } &  1.20  &  1.10  &  2.76  &  5.11 \\
\rowcolor{rgrey} &  0.35  &  0.85  &  2.13  &  3.33 \\
\rowcolor{rgrey} \multirow{-2}{*}{ HR  $\rightleftarrows$  RS } &  0.45  &  1.32  &  3.31  &  7.51 \\
 &  1.00  &  2.11  &  5.26  &  13.36 \\
\multirow{-2}{*}{ HR  $\rightleftarrows$  SI } &  1.00  &  1.64  &  4.10  &  7.45 \\
\rowcolor{rgrey}  &  0.70  &  1.91  &  4.77  &  7.37  \\
\rowcolor{rgrey} \multirow{-2}{*}{ HU  $\rightleftarrows$  RO }  &  0.70  &  3.03  &  7.57  &  15.30\\
 &  0.60  &  1.66  &  4.14  &  6.47 \\
\multirow{-2}{*}{ HU  $\rightleftarrows$  RS } &  0.70  &  2.25  &  5.64  &  12.52 \\
\rowcolor{rgrey}  &  0.60  &  4.14  &  10.35  &  22.49 \\
\rowcolor{rgrey} \multirow{-2}{*}{ HU  $\rightleftarrows$  SK } &  1.30  &  3.24  &  8.09  &  13.72 \\
 &  0.16  &  2.20  &  5.50  &  12.96 \\
 \multirow{-2}{*}{ IT  $\rightleftarrows$  SI } &  0.58  &  2.08  &  5.21  &  13.70 \\
\rowcolor{rgrey} &  3.60  &  2.23  &  5.56  &  9.12 \\
\rowcolor{rgrey} \multirow{-2}{*}{ NO  $\rightleftarrows$  SE } &  3.90  &  1.99  &  4.98  &  9.14 \\
 &  0.60  &  2.45  &  6.13  &  10.11 \\
\multirow{-2}{*}{ PL  $\rightleftarrows$  SK } &  0.50  &  2.99  &  7.47  &  14.22 \\
\rowcolor{rgrey}  &  0.70  &  1.30 &  3.24  &  5.68 \\
\rowcolor{rgrey} \multirow{-2}{*}{ RO  $\rightleftarrows$  RS } &  0.50  &  0.78  &  1.94  &  3.77 \\
    \bottomrule
    \end{tabularx}
    \caption* {Table 3.b: Interconnector capacities cont. (NRL, No Realistic Limit)}
\end{table*}

\end{document}